\input amstex
\documentstyle{amsppt}
\TagsOnRight \NoRunningHeads
\magnification=\magstep1 \vsize 22 true cm \hsize 16 true cm

\topmatter
\title On root-class residuality of generalized free products
\endtitle
\author D.~N.~Azarov and D.~Tieudjo \endauthor
\abstract Root-class residuality of free product of root-class
residual groups is demonstrated. A sufficient condition for
root-class residuality of generalized free product $G$ of groups
$A$ and $B$ amalgamating subgroups $H$ and $K$ through the
isomorphism $\varphi$ is derived. For the particular case when
$A=B$, $H=K$ and $\varphi$ is the identity mapping, it is shown
that group $G$ is root-class residual if and only if $A$ is
root-class residual and subgroup $H$ of $A$ is root-class closed.
These results are extended to generalized free product of
infinitely many groups amalgamating a common subgroup.
\endabstract

\endtopmatter

\document

\centerline {\bf 1. ~ Introduction}

\smallskip

Let $\Cal{K}$ be an abstract class of groups containing at least
one non-identity group. Then $\Cal{K}$ is called a root-class if
the following conditions are satisfied:

1. If $A \in \Cal {K}$ and $B \leqslant A$, then $B \in \Cal {K}$.

2. If $A \in \Cal {K}$ and $B \in \Cal {K}$, then $A \times B \in
\Cal {K}$.

3. If $1 \leqslant C \leqslant B \leqslant A$ is a subnormal
sequence and $A/B, \quad B/C \in \Cal {K}$, then there exists a
normal subgroup $D$ in group $A$ such that $D \leqslant C$ and
$A/D \in \Cal {K}$. See for example $[6]$, p. 428 for details
about this definition.

In this paper, we study root-class residuality of generalized free
products.

We recall that a group $G$ is root-class residual (or $\Cal
{K}$-residual, for a root-class $\Cal {K}$) if, for every
non-identity element $g \in G$, there exists an homomorphism
$\varphi$ from $G$ to some group $G'$ of root-class $\Cal {K}$
such that $g\varphi \ne 1$. Equivalently, $G$ is $\Cal
{K}$-residual if, for every non-identity element $g \in G$, there
exists a normal subgroup $N$ of $G$ such that $G/N \in \Cal {K}$
and $g \notin N$. The most investigated residual properties of
groups are residual finiteness and finite $p$-groups residuality,
(i.e. residuality by the classes of all finite groups and all
finite $p$-groups respectively), and also residuality by the class
of all soluble groups. All these three classes are root-classes.
Therefore results about root-class residuality have sufficiently
enough general character.

In $[6]$ (p. 429) the following result obtained by Gruenberg is
given:

{\sl Free product of root-class residual groups is root-class
residual if and only if every free group is root-class residual.}

The following theorem shown in item 2 asserts that the given above
necessary and sufficient condition is satisfied for any
root-class:

\proclaim {\indent Theorem 1} Every free group is root-class
residual.
\endproclaim

So, Gruenberg's result can be reformulated as follows:

\proclaim {\indent Theorem 2} Free product of root-class residual
groups is root-class residual.
\endproclaim

From theorem 2 and H. Neumann's theorem ([5], p. 122), the
following result is easily established:

\proclaim {\indent Theorem 3} The generalized free product $G$ of
groups $A$ and $B$ amalgamating subgroup $H$ is root-class
residual if groups $A$ and $B$ are root-class residual and there
exists an homomorphism $\varphi$ from $G$ to a group $G'$ of a
root-class such that $\varphi$ is injective on $H$.
\endproclaim

Let's remark that theorem 2 can be considered as a particular case
of theorem 3.

We also see that, if the amalgamated subgroup $H$ is finite, then
the formulated above sufficient condition of root-class
residuality of group $G$ will be as well necessary.

Another restriction permitting to obtain simple criteria of
root-class residuality of generalized free product of groups $A$
and $B$ amalgamating subgroup $H$ is the equality of the free
factors $A$ and $B$.

More precisely, let $G$ be the generalized free product of groups
$A$ and $B$ amalgamating subgroups $H$ and $K$ through the
isomorphism $\varphi$. If $A=B$, $H=K$ and $\varphi$ is the
identity map, we denote group $G$ by $Q = A \underset {H} \to\star
A$. We call $Q$ the generalized free square of group $A$ over
subgroup $H$. Then for such group $Q$ we prove the following
criterium:

\proclaim {\indent Theorem 4} Group $Q = A \underset {H} \to\star
A$ is root-class residual if and only if group $A$ is root-class
residual and subgroup $H$ of $A$ is root-class closed.
\endproclaim

In $[3]$ the above result is obtained for the particular case of
the class of all finite $p$-groups. We recall that subgroup $H$ of
a group $A$ is root-class closed (or $\Cal {K}$-closed, for a
root-class $\Cal {K}$) if, for any element $a$ of $A$ and $a
\notin H$, there exists an homomorphism $\varphi$ from $A$ to a
group of root-class $\Cal {K}$ such that $a\varphi \notin
H\varphi$. This means that, for each $a \in A \setminus H$, there
exists a normal subgroup $N$ of $A$ such that $A/N \in \Cal {K}$
and $a \notin NH$. A $p$-closeness analogue of this definition is
given in $[4]$.

Further, the generalized free product of infinitely many groups
amalgamating subgroup is introduced in $[7]$. Some results on
residual properties of this construction are shown in $[1]$. We
extend theorems 3 and 4 above to generalized free products of
every family $(G_\lambda)_{\lambda \in \Lambda}$ of groups
$G_\lambda$ amalgamating a common subgroup $H$ (theorems 5 and 6).
The set $\Lambda$ can be infinite. Theorems 5 and 6 generalize
some of the results obtained in $[1]$.

\newpage

\centerline {\bf 2. Proofs of theorems 1-4}

\smallskip

Let $\Cal {K}$ be a root-class of groups.

\proclaim {\indent Lemma} Then

1. If a group $G$ has a subnormal sequence with factors belonging
to class $\Cal {K}$, then $G \in \Cal {K}$.

2. If $F \trianglelefteq G,\ G/F \in \Cal {K}$ and $F$ is
 $\Cal {K}$-residual, then group $G$ is also $\Cal {K}$-residual.

3. If $A \trianglelefteq G,\ B \trianglelefteq G,\ G/A \in \Cal
{K}$ and $G/B \in \Cal {K}$, then $G / (A \cap B) \in \Cal {K}$.
\endproclaim

In fact, from the definition of root-class, it follows that
root-class is closed for extensions. So the first property of
lemma is satisfied. The second and third properties are also
easily verified by the definition of root-class.

\vskip 0.3 cm

For the proof of theorem 1, let's remark that every root-class
$\Cal {K}$ contains a non-identity cyclic group (property 1 of the
definition of root-class). If $\Cal {K}$ contains an infinite
cyclic group then, by lemma, $\Cal {K}$ contains any group
possessing subnormal sequence with infinite cyclic factors; and
thus all finitely generated nilpotent torsion-free groups belong
to class $\Cal {K}$. If $\Cal {K}$ contains a finite non-identity
cyclic group, then $\Cal {K}$ contains group of prime order $p$
and consequently, by lemma, $\Cal {K}$ contains all groups
possessing subnormal sequence with factors of order $p$; hence all
finite $p$-groups belong to $\Cal {K}$. So any root-class contains
all finitely generated nilpotent torsion-free groups or all finite
$p$-groups for some prime $p$. Therefore, to end the proof of
theorem 1, let's remind that free groups are residually finitely
generated nilpotent torsion-free groups and also residually finite
$p$-groups (see $[2]$, p. 121 and $[6]$, p. 347).

\vskip 0.3 cm

From the proof of theorem 1 and the Grunberg's theorem formulated
above theorem 2 directly follows.

\vskip 0.3 cm

Let's prove theorem 3.

Let $G$ be the generalized free product of groups $A$ and $B$
amalgamating subgroup $H$ and let groups $A$ and $B$ be $\Cal
{K}$-residual. Suppose there exists an homomorphism $\sigma$ of
$G$ to a group of class $\Cal {K}$, which is one-to-one on $H$.
Let's denote by $N$ the kernel of the homomorphism $\sigma$. Then
$G/N \in \Cal {K}$ and $N \cap H = 1$. By H. Neumann's theorem
([5], p. 122) $N$ is the the free product of a free group $F$ and
some subgroups of group $G$ of the form
$$
g ^ {-1} Ag \cap N, \quad g ^ {-1} Bg \cap N, \tag {$ \ast $}
$$
where $g \in G$. The subgroups of the form $(*)$ are root-class
residual since are groups $A$ and $B$. By theorem 1, free group
$F$ is also root-class residual. Thus $N$ is a free product of
root-class residual groups. Therefore, by theorem 2, $N$ is
root-class residual. Moreover, since $G/N \in \Cal {K}$, by
property 2 of lemma, it follows that group $G$ is root-class
residual. Theorem 3 is proven.

\vskip 0.3 cm

Let's now prove theorem 4.

Let $Q = A \underset {H} \to\ast A$. For any normal subgroup $N$
of group $A$ one can define the generalized free square
$$
Q_N=A/N \underset {HN/N} \to\ast A/N
$$
of group $A/N$ over subgroup $HN/N$ and the homomorphism
$\varepsilon_N: Q \longrightarrow Q_N$, extending the canonical
homomorphism $A \longrightarrow A/N$. It is evident that group
$Q_N$ is an extension of free group with group $A/N$. So, if $A/N$
belongs to root-class $\Cal {K}$ then, by lemma and theorem 1,
$Q_N$ is $\Cal {K}$-residual. Thus, to end the proof, it is enough
to show that $Q$ is residually a group of the form $Q_N$ such that
$A/N \in \Cal {K}$.

Suppose group $A$ is $\Cal {K}$-residual and subgroup $H$ of $A$
is $\Cal {K}$-closed. Let $g \in Q$ such that $g \ne 1$. And let
$g=a_1 \cdots a_s$ be the irreducible form of element $g$. Two
cases arise:

1. $s > 1$. In this case $a_i \in A\setminus H$ for all $i=1,
\ldots, s$. From $\Cal {K}$-closeness of $H$, it follows that, for
every $i=1, \cdots, s$, there exits a normal subgroup $N_i$ of
group $A$ such that $A/N_i \in \Cal {K}$ and $a_i \notin HN_i$.
Let $N=N_1 \cap \cdots \cap N_s$. By lemma, $A/N \in \Cal {K}$
and, it is clear that, for all $i=1, \cdots, s$, $a_i \notin HN$
i.e. $a_iN \notin HN/N$. So, for all $i=1, \cdots, s$, $a_i
\varepsilon_N \notin H \varepsilon_N$. Therefore the form
$$
g \varepsilon_N = a_1 \varepsilon_N \cdots a_s \varepsilon_N
$$
is irreducible and has length $s > 1$.

\noindent Consequently $g \varepsilon_N \ne 1$.

2. $s=1$ i.e. $g \in A $. As group $A$ is $\Cal {K}$-residual,
there exists a normal subgroup $N$ of $A$ such that $A/N \in \Cal
{K}$ and $g \notin N$, i.e. $gN \ne N$. Hence $g \varepsilon_N \ne
1$.

Thus, in any case, for an element $g \ne 1$ in group $A$, there
exists a normal subgroup $N$ such that $A/N \in \Cal {K}$ and the
homomorphism $\varepsilon_N: Q \longrightarrow Q_N$ transforms $g$
to a non identity element. Hence group $Q$ is residually a group
$Q_N$ where $A/N \in \Cal {K}$. Therefore $Q$ is $\Cal
{K}$-residual.

Conversely, suppose group $Q$ is $\Cal {K}$-residual. Evidently
his subgroup $A$ has the same property. Let's prove that $H$ is a
$\Cal {K}$-closed subgroup of group $A$. Let $\gamma$ be an
automorphism of group $Q$ canonically permuting the free factor.
Let $a \in A\setminus H$. Then $a \gamma \ne a$. As $Q$ is $\Cal
{K}$-residual, there exists a normal subgroup $N$ of $Q$ such that
$Q/N \in \Cal {K}$ and $aN \ne a \gamma N$. Let $M=N \cap N
\gamma$. Then
$$
M \gamma = N \gamma \cap N \gamma^2 = N \gamma \cap N =M.
$$
Consequently, in the quotient-group $Q/M$, it is possible to
consider the automorphism $\overline {\gamma}$, induced by
$\gamma$. As $aN \ne a \gamma N$ and $M \leqslant N$, $aM \ne a
\gamma M$. On the other hand, $a \gamma M = (aM) \overline
{\gamma}$. Thus $aM \ne (a M) \overline {\gamma}$. Since $\gamma$
acts identically on $H$ then $\overline\gamma$ also acts
identically on $HM/M$. So and since $aM \ne (aM) \overline\gamma$,
it follows that $aM\notin HM/M$ i.e. $a \varepsilon \notin H
\varepsilon$, where $\varepsilon$ is the canonical homomorphism of
group $Q$ onto $Q/M$. Consequently, $Q/M \in \Cal {K}$ and the
$\Cal {K}$-closeness of subgroup $H$ of group $A$ is demonstrated.

\vskip 0.3 cm

Let's remark, in summary, that the necessary condition for theorem
4 takes place even at more gentle restriction on class $\Cal K$,
namely when $\Cal K$ satisfies only properties 1 and 2 of the
definition of root-class.

\vskip 0.5cm

\centerline {\bf 3. ~ Generalization}

\smallskip

Let $(G_\lambda)_{\lambda \in \Lambda}$ be a family of groups,
where the set $\Lambda$ can be infinite. Let $H_\lambda \leqslant
G_\lambda$, for every $\lambda \in \Lambda$. Suppose also that,
for every $\lambda, \mu \in \Lambda$, there exists an isomorphism
$\varphi_{\lambda \mu}$ : $H_\lambda \longrightarrow H_\mu$ such
that, for all $\lambda, \mu, \nu \in \Lambda$, the following
conditions are satisfied: $\varphi_{\lambda \lambda} =
id_{H_\lambda}$, $\varphi^{-1}_{\lambda \mu} = \varphi_{\mu
\lambda}$, $\varphi_{\lambda \mu} \varphi_{\mu \nu} =
\varphi_{\lambda \nu}$. Let now
$$
G = (G_\lambda\ (\lambda \in \Lambda);\ h \varphi_{\lambda \mu} =
h\ (h \in H_\lambda,\ \lambda, \mu \in \Lambda))
$$
be the group generated by groups $G_\lambda\ (\lambda \in
\Lambda)$ and defined by all the relators of these groups and
moreover by all possible relations: $h \varphi_{\lambda \mu} = h$,
where $h \in H_\lambda,\ \lambda, \mu \in \Lambda$. It is evident
every $G_\lambda$ can be canonically embedded in group $G$ and if
we consider $G_\lambda \leqslant G$ then, for all different
$\lambda, \mu \in \Lambda$,
$$
G_\lambda \cap G_\mu = H_\lambda = H_\mu .
$$
Let's denote by $H$ the subgroup of group $G$, equal to the common
subgroups $H_\lambda$. Then $G$ is the generalized free product of
the family $(G_\lambda)_{\lambda \in \Lambda}$ of groups
$G_\lambda\ (\lambda \in \Lambda)$ amalgamating subgroup $H$. We
will consider, as well, that $G_\lambda \leqslant G$, for all
$\lambda \in \Lambda$. See $[1]$ or $[7]$ for details about
generalized free product of a family of groups.

\proclaim {\indent Theorem 5} The generalized free product $G$ of
the family $(G_\lambda)_{\lambda \in \Lambda}$ of groups
$G_\lambda$ amalgamating subgroup $H$ is root-class residual if
every group $G_\lambda$ is root-class residual and there exists an
homomorphism $\varphi$ from $G$ to a group $G'$ of a root-class
such that $\varphi$ is injective on $H$.
\endproclaim

\demo{Proof} The proof is the same as that of theorem 3.

In fact, let groups $G_\lambda$ be $\Cal {K}$-residual, for all
$\lambda \in \Lambda$. Suppose there exists an homomorphism
$\sigma$ of $G$ to a group of class $\Cal {K}$, which is
one-to-one on $H$ and let $N = ker \sigma$. Then $G/N \in \Cal
{K}$ and $N \cap H = 1$. But $N$ is the the free product of a free
group $F$ and some subgroups of group $G$ of the form
$$
g ^ {-1} G_\lambda g \cap N, \tag {$ \ast \ast $}
$$
(where $g \in G$ and $\lambda \in \Lambda$) which are root-class
residual. Since $F$ is also root-class residual by theorem 1, $N$
is a free product of root-class residual groups. Thus, by theorem
2, $N$ is root-class residual. Moreover, since $G/N \in \Cal {K}$,
by property 2 of lemma, it follows that group $G$ is root-class
residual and the theorem is proven.
\enddemo

\bigskip

Suppose now that, for all $\lambda \in \Lambda$, $G_\lambda = A$
and denote $P$, the generalized free power of group $A$ over
subgroup $H$, by $A \underset {H} \to{\star} \cdots \underset {H}
\to{\star} A$. For such group $P$ we have the following criterium:

\proclaim {\indent Theorem 6} Group $P=A \underset {H} \to{\star}
\cdots \underset {H} \to{\star} A$ is root-class residual if and
only if group $A$ is root-class residual and subgroup $H$ of $A$
is root-class closed.
\endproclaim

The proof is similar to that of theorem 4.

\bigskip

Ivanovo State University, Ermaka str. 37, 153025 Ivanovo, Russia.

\smallskip

University of Ngaoundere, P.~O.~BOX 454, Ngaoundere, Cameroon.

{\it E-mail:} tieudjo\@yahoo.com

\end
\vskip 0.5 cm


\Refs

\ref \no 1 \by D. Doniz \paper Residual properties of free
products of infinitely many nilpotent groups amalgamating cycles
\jour J. Algebra \vol 179 \yr 1996 \pages 930--935
\endref

\ref \no 2 \by M. I. Kargapolov, I. I. Merzliakov \book Elements
of group theory \lang Russian \publ M., Naouka  \yr 1972
\endref

\ref \no 3 \by G. Kim and J. McCarron \paper On amalgamated free
products of residually $p$-finite groups \jour J. Algebra, \vol
162 \yr 1993 \pages 1--11 \endref

\ref \no 4 \by G. Kim and C. Y. Tang \paper On generalized free
products of residually finite $p$-groups \jour J. Algebra, \vol
201 \yr 1998 \pages 317--327 \endref

\ref \no 5 \by  R. Lyndon, P. Schupp \book Combinatorial group
theory \lang Russian \publ M., Mir \yr 1980 \endref

\ref \no 6 \by W. Magnus, A. Karrass and D. Solitar \book
Combinatorial group theory \lang Russian \publ M., Naouka \yr 1974
\endref

\ref \no 7 \by M. Shirvani \paper A converse to a residual
finiteness theorem of G. Baumslag \jour Proc. Amer. Math. Soc.
\vol 104 \yr 1988 \issue 3 \pages 703--706 \endref


\endRefs
